\documentclass{ifacconf}

\usepackage{graphicx}      
\usepackage{natbib}        
\usepackage{amsmath}
\usepackage{amssymb}
\usepackage{tikz}
\usepackage{pgfplots}
\usepackage{float}
\usepackage{booktabs}
\usepackage[scientific-notation=true]{siunitx}
\usepackage[T1]{fontenc}
\newcommand{\R}{\mathbb{R}}
\newcommand{\bx}{\mathbf{x}}
\newcommand{\by}{\mathbf{y}}

\newcommand{\bff}{\mathbf{f}}
\newcommand{\bg}{\mathbf{g}}
\newcommand{\bu}{\mathbf{u}}

\tikzset{every picture/.style={/utils/exec={\sffamily}}}

\usetikzlibrary{pgfplots.groupplots}
\pgfplotsset{compat=1.14}
\usepgfplotslibrary{colorbrewer}
\definecolor{brinkpink}{rgb}{0.98, 0.38, 0.5}

\definecolor{royalblue}{rgb}{0.08, 0.38, 0.95}

\definecolor{graph_green}{HTML}{009e8c}
\definecolor{graph_petrol}{HTML}{133840}
\definecolor{graph_red}{HTML}{F22738}
\definecolor{graph_grey}{HTML}{7a8581}
\definecolor{graph_purple}{HTML}{8254F6}
\usetikzlibrary{decorations.markings}
\usetikzlibrary{arrows}
\usetikzlibrary{shapes, arrows.meta, positioning}
\tikzset{cross/.style={cross out, draw, 
         minimum size=2*(#1-\pgflinewidth), 
         inner sep=0pt, outer sep=0pt}}
\tikzset{
  font={\fontsize{8pt}{12}\selectfont}}
 \pgfplotsset{
    legend style = {font = \scriptsize},
  }

\begin{document}
\begin{frontmatter}

\title{Data-driven initialization of deep learning solvers for Hamilton-Jacobi-Bellman PDEs\thanksref{footnoteinfo}} 

\thanks[footnoteinfo]{DK was supported
	by the UK Engineering and Physical Sciences Research Council (EPSRC) grants
	EP/V04771X/1, EP/T024429/1, and EP/V025899/1.}

\author[First]{A. Borovykh}
\author[Second]{D. Kalise}
\author[Third]{A. Laignelet} 
\author[Third]{P. Parpas}

\address[First]{Warwick Business School, University of Warwick, Coventry CV4 7AL, UK (e-mail: anastasia.borovykh@wbs.ac.uk)}
\address[Second]{Department of Mathematics, Imperial College London, 
	London SW7 2AZ, UK (e-mail: d.kalise-balza@imperial.ac.uk)}
\address[Third]{Department of Computing, Imperial College London, 
   London SW7 2AZ, UK (e-mail: \{alexis.laignelet18, panos.parpas\}@imperial.ac.uk)}

\begin{abstract}                
    A deep learning approach for the approximation of the Hamilton-Jacobi-Bellman partial differential equation (HJB PDE) associated to the Nonlinear Quadratic Regulator (NLQR) problem. A state-dependent Riccati equation control law is first used to generate a gradient-augmented synthetic dataset for supervised learning. The resulting model becomes a warm start for the minimization of a loss function based on the residual of the HJB PDE. The combination of supervised learning and residual minimization avoids spurious solutions and mitigate the data inefficiency of a supervised learning-only approach. Numerical tests validate the different advantages of the proposed methodology.

\end{abstract}


\end{frontmatter}

\section{Introduction}


Over the last years there has been remarkable progress on the solution of high-dimensional nonlinear optimal control problems, see e.g. \cite{dolgov2021tensor,onken2021,walter2021,meng2022}, partly fuelled by the effectiveness of deep learning methods for solving nonlinear PDEs including HJB equations (see e.g. \cite{raissi2017physics, sirignano2018dgm, han2018solving, guler2019towards}). In its general form, \emph{unsupervised} deep learning-based PDE solvers aim at minimizing a loss function in which the interior and boundary constraints from the PDE are incorporated (see e.g. eq. \eqref{eq:residual_loss}).
Certain challenges can still arise; the works of e.g. \cite{van2020optimally, wang2021understanding} have shown that the accuracy of the solution obtained by the neural network heavily relies on the weighting between the different terms (interior, boundary) in the loss function. As also noted in \cite{van2020optimally} a long training time can be required and the number of domain points one has to sample can become large. Furthermore, if the PDE solutions are non-unique, using solely the residual can lead to an incorrect solution; we show this for the linear HJB equation in Fig.~\ref{fig:linear_solutions}. As an alternative, \emph{supervised} learning methods can be employed (see e.g. \cite{azmi2020optimal, liu2019pricing, kang21}) where an optimal control solver is used to generate the solution dataset and a machine learning model is then trained on this set of data in a supervised manner. Alternatively, in \cite{albi2021gradient} the solutions of the state-dependent Riccati equation (SDRE) (see e.g. \cite{banks2007nonlinear}) are used as the initial dataset. A downside of such supervised settings is the propagation of errors from the initial solver into the model learned by the neural network (see Fig.~\ref{fig:non_linear_discrepancy} where the discrepancy between the residual loss and the SDRE-generated solution is highlighted). 

In this work we focus on the synthesis of feedback laws for optimal nonlinear stabilization problems, that is, by solving the HJB PDE associated to the NLQR problem. We propose a neural network-based approximation for the value function where the training is done using a two-step learning approach: in the first step (i) the parameters of the neural network are initialized on the state-value dataset obtained from solving the SDREs, and in the second step (ii) the residual coming from the PDE constraints is minimized. The benefits of this approach include: 
\begin{enumerate}
    \item by first pre-training on the data, our method ensures the convergence to the correct solution; this is particularly relevant in a setting where the solution to the HJB PDE is non-unique (see e.g. Figure \ref{fig:invalid_vs_valid_residual}),
    \item by using the residual-minimization in the second learning step our method allows to mitigate the errors in the SDRE solution and allows to reach higher-accuracy solutions than data-only (see e.g. Figure \ref{fig:mse_data_residual}),
    \item the two-step learning approach decreases the amount of data points required to reach a satisfactory accuracy compared to training on the dataset only (see e.g. Figure \ref{fig:mse_data_residual}).
\end{enumerate}


\section{Methodology}\label{sec:methodology}
\subsection{Infinite horizon nonlinear optimal stabilization}
Let $\by(t)\in\mathbb{R}^n$  and $\bu(t)\in\mathbb{R}^m$ be the state and control signal, respectively, of a system represented in control-affine form
\begin{align}\label{eq:dynamics_y}
	\frac{d}{dt}\mathbf{y}(t)= \mathbf{f}(\by(t))+\mathbf{g}(\by(t))\bu(t), \;\; \by(0) = \bx\,.
\end{align}
We consider the following infinite horizon optimal control problem:
\begin{align}\label{eq:objective_u}
	\min_{\bu(\cdot)\in\mathcal{U}}J(\bu;\bx):=\frac12\int_{0}^{\infty}\by^{\top}Q\by+\bu^{\top}R\bu \,dt\,,
\end{align}
subject to \eqref{eq:dynamics_y}, where the dynamics $\bff:\mathbb{R}^n\rightarrow\mathbb{R}^n$ and $\bg:\mathbb{R}^n\rightarrow\mathbb{R}^{n\times m}$ are assumed to be continuously differentiable, and the running cost matrices are given by $Q\in\mathbb{R}^{n\times n}, Q\succeq 0\,,$ and  $R\in\mathbb{R}^{m\times m}, R\succ 0$. We set $\mathcal{U}\equiv L^2([0,\infty);\R^m)$.

The formulation above corresponds to the NLQR problem, which arises in the synthesis of optimal feedback laws for the stabilization of nonlinear dynamics. The solution of this problem follows standard dynamic programming argument, see e.g., \cite{bardi1997optimal}. Defining the value function of the problem
\begin{align}\label{eq:optimal_cost}
	V(\mathbf{x})=\inf_{\mathbf{u}(\cdot)} J(\mathbf{u};\mathbf{x})\,,
\end{align}
the optimal control is expressed in feedback form as
\begin{align}\label{eq:opt_control}
	\mathbf{u}^*(\mathbf{x})=-R^{-1}\mathbf{g}(\mathbf{x})^{\top}\nabla V(\mathbf{x})\,,
\end{align}
where $V(\bx)$ satisfies the Hamilton-Jacobi-Bellman PDE
\begin{align}\label{eq:hjb_pde}
	&\begin{aligned}
		& \mathcal{N}(\mathbf{x},V):=-\frac{1}{2} \nabla V(\mathbf{x})^{\top}\mathbf{g}(\mathbf{x})R^{-1}\mathbf{g}^{\top}(\mathbf{x})\nabla V(\mathbf{x})\\
		& \qquad + \nabla V(\mathbf{x})^{\top}\mathbf{f}(\mathbf{x})+\frac12\mathbf{x}^{\top}Q\mathbf{x}=0\,,
	\end{aligned}
\end{align}
in addition to the boundary condition $V(\mathbf{0})=0$.
\subsection{Neural network model}
In this work, we are interested in approximating the value function using deep neural networks. For this, we consider a fully connected feedforward neural network given by
\begin{align}
    \hat V (\mathbf{x};\theta) = f_L\circ f_{L-1}\circ\cdots\circ f_1(\mathbf{x})
\end{align}
where $\mathbf{x}$ is the model input, $L$ is the number of layers and $f_i(\mathbf{x})=\sigma(W^i\mathbf{x}+\mathbf{b}^i)$, where $W^i$ and $b^i$ are matrices and vectors of prescribed dimensions, and $\sigma$ is a nonlinear activation function applied component-wise. The unknown parameters of the network $\{W^i,b^i\}_{i=1}^L$ are grouped in the vector $\theta$. The neural network is trained upon a set of collocation points $\{\bx^i\}_{i=1}^N$ sampled uniformly over a domain of interest. The training of the neural network in based on the two-step learning procedure presented in Section \ref{sec:training}.

\subsection{Synthetic data generation for supervised learning}
In our two-step learning procedure $\hat V (\mathbf{x};\theta)$ is first trained following a supervised learning gradient-augmented approach, that is, upon a collection of values \[\{\bx^i,V(\bx^i),\nabla V(\bx^i)\}_{i=1}^{N_1}\,,\] approximating pointwise valuations of the true value function and its gradient. This dataset is generated synthetically by resorting to the link between the value function of the control problem and optimality conditions.  One possibility to generate the dataset in a finite horizon setting is to use Pontryagin's Maximum Principle to obtain first-order optimality conditions for \eqref{eq:dynamics_y}-\eqref{eq:objective_u} in the form of a two-point boundary value problem, separately for each initial condition $\bx^i$. However, in the infinite horizon optimal stabilization problem setting studied in this work, such an interpretation is not readily available, as the two-point boundary value problem has an asymptotic terminal condition for $t\to\infty$. We therefore rely on the use of a state-dependent Riccati equation approach \cite{banks2007nonlinear} which provide a stabilizing feedback control through a sequential solution of the algebraic Riccati equations (ARE) along a trajectory. 

\subsubsection{State-Dependent Riccati Equation}
For nonlinear dynamics, using the state-dependent Riccati equation (SDRE) approach is an effective method to generate a suboptimal, yet stabilizing solution which can be interpreted as an approximation of the HJB PDE \eqref{eq:hjb_pde}   (see e.g. \cite{jones2020solution,albi2021gradient}). This approach is based on an extended linearization of the dynamics, expressed in semilinear form $f(\bx)=A(\bx)\bx$. Without loss of generality let $V(\bx)=\frac{1}{2}\mathbf{x}^{\top}P(\mathbf{x})\mathbf{x}$. Note that $\nabla V(\mathbf{x})=P(\bx)\mathbf{x} + \phi(\mathbf{x})$ where  $\phi\in\mathbb{R}^n$ is a `discrepancy' term arising from the differentiation of the non-linearity in $P(\mathbf{x})$ and its entries are given by
\begin{align*}
    &\phi(\mathbf{x})_k = \sum_{i,j=1}^nx_ix_j\nabla_{x_k}P_{ij}(\mathbf{x})\,,\quad k=1,\ldots,n\,.
\end{align*}

In the SDRE approach, we approximate $\nabla V(\mathbf{x})\approx P(\mathbf{x})\bx$ so that the HJB PDE for the NLQR simplifies to 
\begin{align}\label{eq:sdre}
    \begin{aligned}
    & -\frac{1}{2}P(\mathbf{x})g(\mathbf{x})R^{-1}g(\mathbf{x})^{\top}P(\mathbf{x})+P(\mathbf{x})A(\mathbf{x}) \\
    & \qquad +A(\mathbf{x})^{\top}P(\mathbf{x})+Q=0,
    \end{aligned}
\end{align}
For a fixed $\mathbf{x}$, the above equation amounts to solving an \emph{ARE}. This suggests an implementation in a model predictive control fashion. Specifically, given a current state $\mathbf{\bar x}$ we solve \eqref{eq:sdre} for $P(\mathbf{\bar x})$ by fixing every operator. Then using \eqref{eq:opt_control} gives 
$\mathbf{u}(\mathbf{x})=-g^{\top}(\mathbf{\bar x})P(\mathbf{\bar x})\mathbf{x}$ and consequently using this in the dynamics of \eqref{eq:dynamics_y} gives the next value of the state. Similarly as in the ARE setting, along a trajectory we generate data for $V(\bx)=\frac12 \bx^{\top}P(\bx)\bx$ and $\nabla V(\bx)=P(\bx)\bx$. However, as long as $\phi(\bx)\neq 0$, this is a suboptimal approximation to the solution of the HJB PDE. Through supervised learning, we train a neural network which is subsequently used as a warm start of the training of a neural network which will directly approximate the solution of the HJB PDE.

\subsection{The training methodology}\label{sec:training}


We propose a two-step learning approach. A first neural network is trained using supervised learning with a synthetic dataset generated from the SDRE approach. In a second step, a neural network, initialized with the parameters of the first step, is trained solely for minimizing the residual of the HJB PDE \eqref{eq:hjb_pde}.  The idea behind such a solver is that training first on the synthetic dataset allows the residual minimization to start with a set of weights closer to the basin of attraction of the correct solution. The residual minimisation is then necessary to reduce the approximation errors induced by the data generation from the SDRE approach. 


\subsubsection{Data-driven initialization}
The SDRE approach generates a gradient-augmented dataset $\{\bx^i,V(\bx^i),\nabla V(\bx^i)\}_{i=1}^{N_1}$. Note that the gradient information is generated at virtually no extra computational cost. Consider the following loss function,
\begin{align}\label{eq:data_loss}
	\begin{aligned}
    \mathcal{L}^{\textnormal{dat}}(\theta) = &\lambda_1\frac{1}{N_1} \sum_{i=1}^{N_1}||V(\bx^i)-\hat V(\bx^i;\theta)||_2^2 \\
    &+ \lambda_2\frac{1}{N_1} \sum_{i=1}^{N_1}||\nabla V(\bx^i)-\nabla \hat V(\bx^i;\theta)||_2^2,
    \end{aligned}
\end{align}
where $N_1$ is the number of data points generated using the SDRE solver and $\lambda_1$ and $\lambda_2$ are the weights associated to the different loss function terms.
Our first step then amounts to initialising the neural network parameters by finding, 
\begin{align}
    \theta^{\textnormal{dat}} = \arg\min_{\theta} \mathcal{L}^{\textnormal{dat}}(\theta). 
\end{align}
We remark that $\nabla \hat V(\bx;\theta)$ is obtained by differentiating the neural network with respect to the input $\mathbf{x}$. 

\subsection{Residual loss function} 
Based on \eqref{eq:hjb_pde}, the \emph{residual} loss function is given by,
\begin{align}\label{eq:residual_loss}
    \mathcal{L}^{\textnormal{res}}(\theta) = \frac{1}{N_2}\sum_{i=1}^{N_2}||\mathcal{N}(\mathbf{x}^i,\hat V(\bx^i;\theta))||_2^2,
\end{align}
where $N_2$ is the number of collocation points $\mathbf{x}^1,...,\mathbf{x}^{N_2}$ which we sample on the interior of the domain. 
As before, the gradient terms in this loss function are obtained by differentiating the neural network output $\hat V$ with respect to state variables. 
After obtaining $\theta^{dat}$ from minimizing \eqref{eq:data_loss}, these parameters are used as the initialization of the training of a second neural network based on the residual minimization:
\begin{align}
    \theta^{\textnormal{res\_valid}} = \arg\min_{\theta} \mathcal{L}^{\textnormal{res}}(\theta)\,.
\end{align}


\subsection{Benefits of the two-step learning approach}\label{benefits}
\begin{figure}
    \centering
\begin{tikzpicture}
\begin{groupplot}[
    group style={
        group name=my plots,
        group size=2 by 4,
        ylabels at=edge left,
        horizontal sep=1.8cm
    },
    width=4.3cm,
    height=4.3cm,
]

\nextgroupplot[
    xlabel=$x_1$,
    ylabel=$x_2$,
    zlabel=$z$,
    view={30}{40},
    grid=both,
    mesh/ordering=x varies,
    mesh/rows=20,
    zmin=0,
    zmax=1.]
     \addplot3[
     opacity=0.8,
     colormap/viridis,
     surf,
     ] table[x index=0,y index=1, z index=2,col sep=comma] {data/lqr_true_solutions.dat};
         
\nextgroupplot[
    xlabel=$x_1$,
    ylabel=$x_2$,
    zlabel=$z$,
    view={30}{40},
    grid=both,
    mesh/ordering=x varies,
    mesh/rows=20,
    zmin=-0.5,
    zmax=0.5]
     \addplot3[
     opacity=0.8,
     colormap/viridis,
     surf,
     ] table[x index=0,y index=1, z index=3,col sep=comma] {data/lqr_true_solutions.dat};

\nextgroupplot[
    xlabel=$x_1$,
    ylabel=$x_2$,
    zlabel=$z$,
    view={30}{40},
    grid=both,
    mesh/ordering=x varies,
    mesh/rows=20,
    zmin=-0.5,
    zmax=0.5,
    ]
     \addplot3[
     opacity=0.8,
     colormap/viridis,
     surf,
     ] table[x index=0,y index=1, z index=5,col sep=comma] {data/lqr_true_solutions.dat};
     
\nextgroupplot[
    xlabel=$x_1$,
    ylabel=$x_2$,
    zlabel=$z$,
    view={30}{40},
    grid=both,
    mesh/ordering=x varies,
    mesh/rows=20,
    zmin=-0.5,
    zmax=0.5,
    ]
     \addplot3[
     opacity=0.8,
     colormap/viridis,
     surf,
     ] table[x index=0,y index=1, z index=4,col sep=comma] {data/lqr_true_solutions.dat};

\end{groupplot}

\end{tikzpicture}
    \caption{Possible solutions of the Hamilton-Jacobi equation in the LQR case with $A = B = Q = I$, and $R = 0.2I$. }
    \label{fig:linear_solutions}
\end{figure}
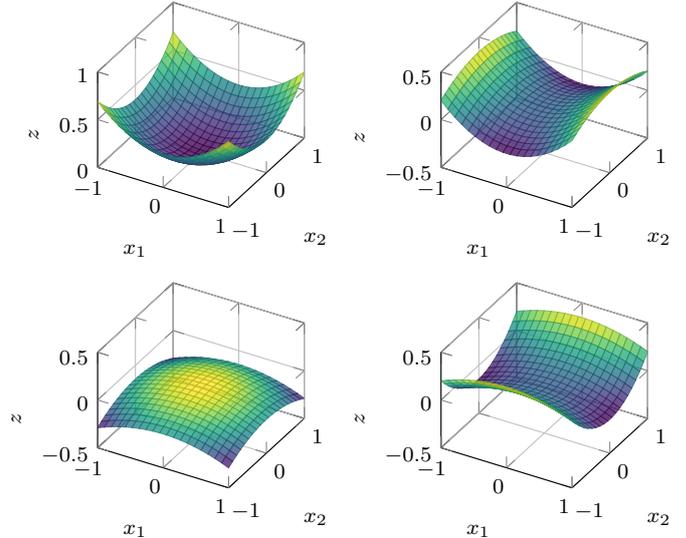

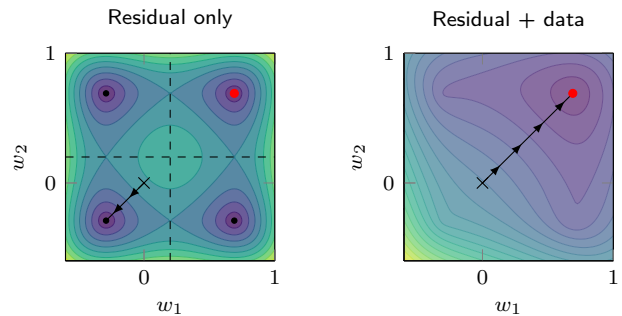
\begin{figure}
 \centering
 \begin{minipage}[b]{0.49\linewidth}
 \begin{tikzpicture}
      \begin{axis}[
        xlabel=$w_1$,
        ylabel=$w_2$,
        zlabel=$z$,
        title=Residual only,
        view={0}{90},
        width=\textwidth, height=\textwidth,
        mesh/ordering=x varies,
        mesh/rows=40,
        zmin=0.0,
        zmax=1.,
      ]
    \addplot[opacity=1] graphics[xmin=-0.6,xmax=1,ymin=-0.6,ymax=1] {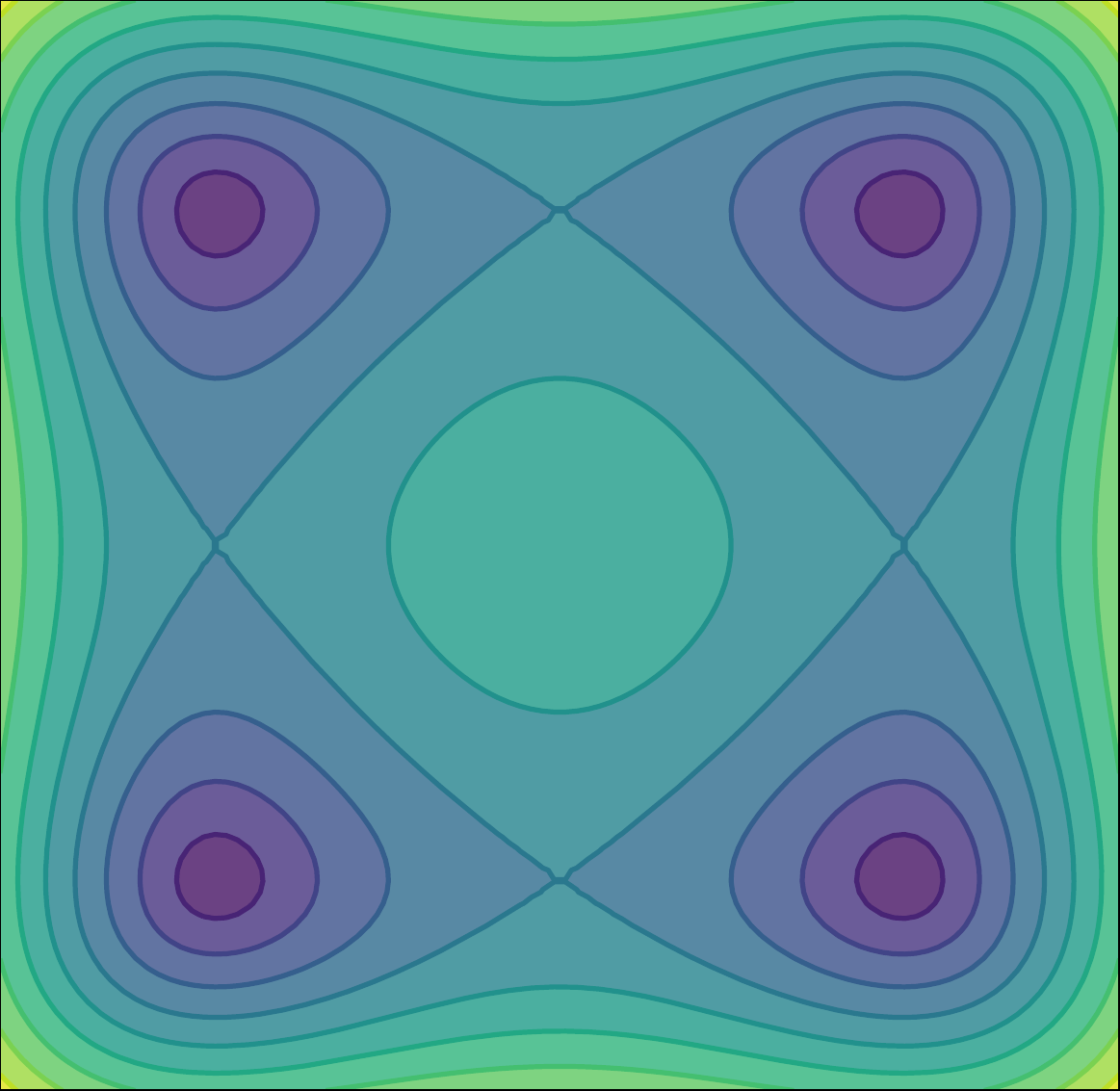};
    \addplot[only marks, mark=*, mark size = 1, black]
        coordinates{ 
          (-0.29,-0.29)
          (0.69,-0.29)
          (-0.29,0.69)
        };
    \addplot[only marks, mark=*, mark size = 1.5, red]
        coordinates{ 
          (0.69,0.69)
        };
    \addplot[only marks, mark=x, mark size = 3, black]
        coordinates{ 
          (0,0)
        };
    \path[ draw,%
      decoration={%
        markings,%
        mark=between positions 0.4 and 0.9 step 3mm with \arrow{latex},%
      },%
      postaction=decorate, black] (0,0) to (-0.29,-0.29);
   \path[ draw, black, dashed] (0.2,-0.6) to (0.2,1);  
   \path[ draw, black, dashed] (-0.6, 0.2) to (1,0.2);  
      \end{axis}
      \end{tikzpicture}
  \end{minipage}
 \begin{minipage}[b]{0.49\linewidth}
 \begin{tikzpicture}
      \begin{axis}[
        xlabel=$w_1$,
        ylabel=$w_2$,
        zlabel=$z$,
        view={0}{90},
        title=Residual + data,
        width=\textwidth, height=\textwidth,
        mesh/ordering=x varies,
        mesh/rows=40,
        zmin=0.0,
        zmax=0.4
      ]
    \addplot[opacity=0.8] graphics[xmin=-0.6,xmax=1,ymin=-0.6,ymax=1] {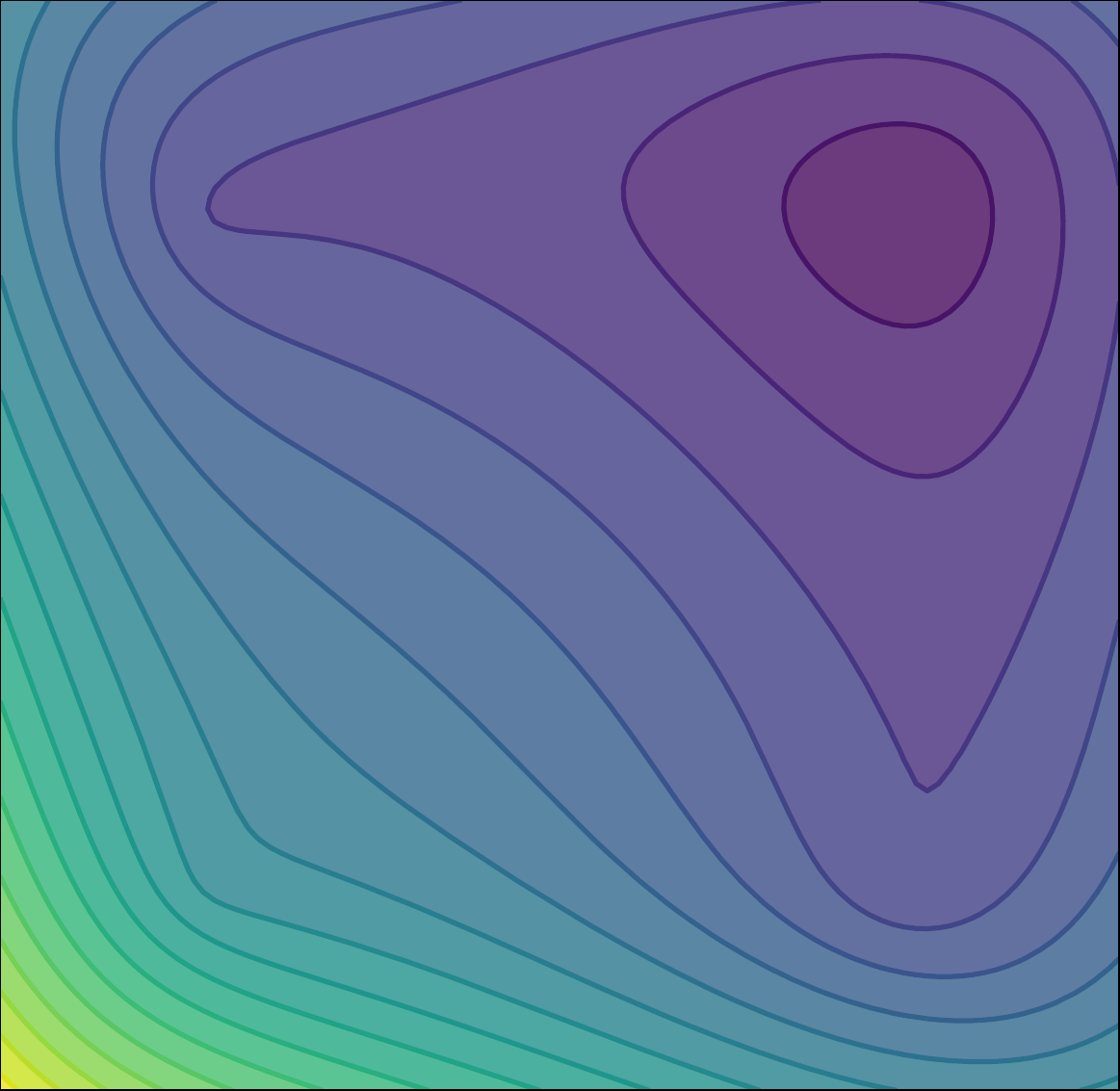};
    \addplot[only marks, mark=*, mark size = 1.5, red]
        coordinates{ 
          (0.69,0.69)
        };
        \addplot[only marks, mark=x, mark size = 3, black]
        coordinates{ 
          (0,0)
        };
        \path[ draw,%
          decoration={%
            markings,%
            mark=between positions 0.2 and 0.95 step 4mm with \arrow{latex},
          },%
          postaction=decorate, black] (0,0) to (0.69,0.69);
      \end{axis}
      \end{tikzpicture}
  \end{minipage}
   \caption{Loss functions in the case of residual-only and a linear combination of the residual and the data term. The global minimum for the valid value function is represented by a red dot. A loss function combining the data and residual has a unique solution, while a residual-only approach inherits the multiple, but wrong, solutions of the original HJB PDE.}
  \label{fig:lqr_loss_function}  
\end{figure}

We illustrate, using a toy example, the benefits of warm start for residual minimization. Consider the linear quadratic control problem with $A=I$, $B=I$, $Q=I$ and $R=0.2I$ with $n=2$. The solution $P$ of the corresponding ARE is a diagonal matrix with elements $\alpha_1 \approx \pm 0.69$ and/or $\alpha_2 \approx\pm -0.29$, see Figure \ref{fig:linear_solutions} for the different solutions $V(\bx)$. By the positive-definiteness assumption the only valid solution is $P=\alpha_1 I$. However, if a neural network was trained on the residual only, four global minima exist (the solutions for $P$ with $P$'s elements different combinations of the roots) and convergence towards the value function will depend on a correct initialization\footnote{Note that in dimension $n$, the number of possible solutions in this example grows to $2^n$.}.  
Using Xavier initialization \cite{glorot2010understanding} the initial weights are centred around zero and hence gradient descent is most likely to converge to the minimum at $(-0.29, -0.29)$ (see the left-hand side of Fig~\ref{fig:lqr_loss_function}). 
To ensure we converge to the valid solution, we can enrich the loss function with a data-driven term, similar to $\mathcal{L}^{\textnormal{dat}}(\theta)$; this term guides convergence to the correct solution at $(0.69, 0.69)$ (see the right-hand side of Fig.~\ref{fig:lqr_loss_function}).

This toy example motivates our combined approach of data and residual loss minimization. Using a residual-only approach inherits the multiplicity of (wrong) solutions of the HJB PDE associated to the NLQR problem. The downside of using both data and residual terms in \emph{a single} loss function is the problem of choosing the right weights between the residual and data term (\cite{van2020optimally}). 
Even more crucially, in non-linear settings the  approximate value function obtained from the SDRE dataset can significantly differ from the value function that satisfies the residual (see also Figure \ref{fig:non_linear_discrepancy}); in such a setting, minimizing the data term and the residual term \emph{at the same time} is not possible. To ensure our framework is robust and efficient in nonlinear settings, we thus propose the two-step learning approach. 

\section{Experiments}\label{sec:numerics}
\subsection{Implementation details}
\vspace{-0.2cm}
\subsubsection{Neural network configuration and loss function}
The architecture we consider is a fully connected neural network with the sigmoid activation function $\sigma(x)=1/(1+e^{-x})$ with 3 hidden layers. For the two-dimensional problem we set 20 neurons per layer while for the 10-dimensional case we use 50 neurons per layer. The neural network is trained using the Adam optimizer \cite{kingma2014adam}.

For the 2D examples, 20 SDRE data points are sampled for supervised learning and 50 collocation points are used to minimize the HJB PDE residual. The training procedure on the dataset is the following: we use $10^3$ iterations with a learning rate of $0.01$ followed by $10^{3}$ iterations with a learning rate of $10^{-3}$. When training simultaneously on the data and gradient we set $\lambda_1,\lambda_2=1$. 
For the residual minimization, the training is first done for $2\times 10^3$ iterations with a learning rate of $0.01$, followed by $4\times 10^{3}$ iterations with a learning rate of $10^{-3}$.

For the Cucker-Smale model, 500  data points are sampled, and $5\times 10^3$ points are used for training on the residual. The pre-training uses $2\times 10^3$ iterations with  a learning rate of $0.01$ followed by $4\times 10^3$ iterations with $10^{-3}$. The residual minimization is done using $10^4$ iterations with  a learning rate of $10^{-3}$ followed by $2\times 10^4$ iterations with $10^{4}$.

\subsubsection{Data generation}
For the 2D examples, the ARE/SDRE is solved analytically to obtain a closed-form solution for $P(\mathbf{x})$ which is then used in the expressions for $V(\mathbf{x})$ and $\nabla V(\mathbf{x})$. 
In the case of the Cucker Smale model, the SDRE is solved numerically using the SciPy Linear Algebra package and $\nabla V(\mathbf{x})=P(\mathbf{x})\mathbf{x}$ is used for the derivative. The collocation points are uniformly sampled in the space.

\subsubsection{Evaluation}
The evaluation of the accuracy of the final solution is done by computing the mean squared error (MSE) between the neural network on newly generated interior domain points and the exact solution. For the 2D tests, the exact solution is obtained either through a closed-form solution, whenever available, or using the numerical solution of a semi-Lagrangian scheme over a fine grid \cite{alla2015}. The final MSE's are computed as the median over 10 runs of the optimization. 

\subsection{A 2D linear example}
\begin{figure}
    \centering
\begin{tikzpicture}
\begin{groupplot}[
    group style={
        group name=my plots,
        group size=2 by 1,
        ylabels at=edge left
    },
    width=4.75cm,
    height=5cm,
    tickpos=left,
    ylabel shift = -5 pt,
    ymode=log,
    xlabel=Iterations
]

\nextgroupplot[title={2D linear},ylabel=MSE]
    \addplot[color=graph_green, thick]
    table[x index=0,y index=1,col sep=comma] {data/linear_all_methods.dat};            
    \addplot[color=graph_red, thick]
    table[x index=0,y index=2,col sep=comma] {data/linear_all_methods.dat};    
    \addplot[color=graph_petrol, thick]
    table[x index=0,y index=4,col sep=comma] {data/linear_all_methods.dat};
            
\nextgroupplot[title={2D nonlinear},legend style={at={(0.8,-0.3)}, anchor=north east,legend columns=3}]
    \addplot[color=graph_green, thick]
    table[x index=0,y index=1,col sep=comma] {data/non_linear_all_methods.dat};            
    \addplot[color=graph_red, thick]
    table[x index=0,y index=2,col sep=comma] {data/non_linear_all_methods.dat};    
    \addplot[color=graph_petrol, thick]
    table[x index=0,y index=4,col sep=comma] {data/non_linear_all_methods.dat};
    
    \addlegendentry{Data}
    \addlegendentry{Data/gradient}
    \addlegendentry{Residual}

\end{groupplot}
\end{tikzpicture}
    \caption{Mean squared error of different training methods for 2D linear and 2D non linear case (computed over $100\times 100$ points in the domain). 
    }
    \label{fig:mse_training_methods}
\end{figure}
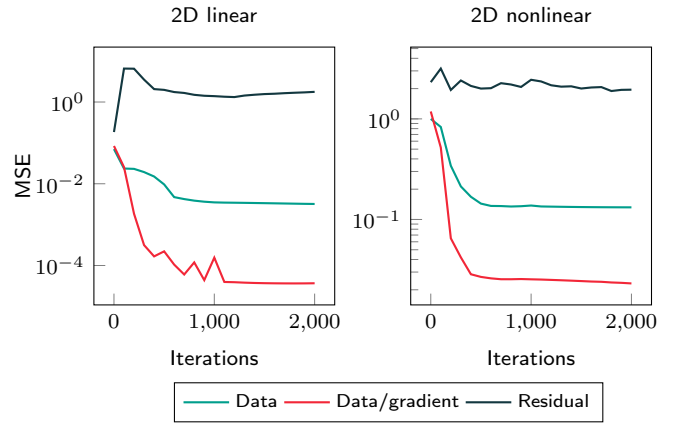

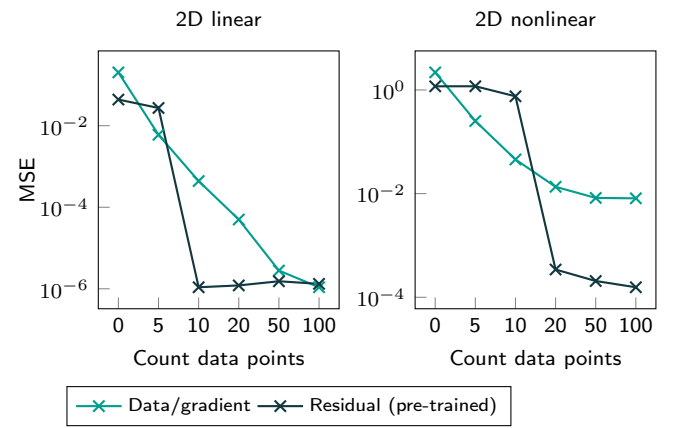
\begin{figure}
    \centering
\begin{tikzpicture}
\begin{groupplot}[
    group style={
        group name=my plots,
        group size=2 by 1,
        ylabels at=edge left
    },
    width=4.75cm,
    height=5cm,
    tickpos=left,
    ylabel shift = -5 pt,
    ymode=log,
    xtick={0,1,2,3,4,5},
    xticklabels={0, 5, 10, 20, 50, 100},
    xlabel=Count data points
]

\nextgroupplot[title={2D linear},ylabel=MSE]
    \addplot[mark=x, mark size=3pt, color=graph_green, thick]
    table[x index=0,y index=2,col sep=comma] {data/linear_results.dat};            
    \addplot[mark=x, mark size=3pt, color=graph_petrol, thick]
    table[x index=0,y index=3,col sep=comma] {data/linear_results.dat};      

\nextgroupplot[title={2D nonlinear},legend style={at={(0.4,-0.3)}, anchor=north east,legend columns=2}]
    \addplot[mark=x, mark size=3pt, color=graph_green, thick]
    table[x index=0,y index=2,col sep=comma] {data/non_linear_results.dat};            
    \addplot[mark=x, mark size=3pt, color=graph_petrol, thick]
    table[x index=0,y index=3,col sep=comma] {data/non_linear_results.dat};      
    
    \addlegendentry{Data/gradient}
    \addlegendentry{Residual (pre-trained)}

\end{groupplot}
\end{tikzpicture}

    \caption{Mean squared error between training on dataset and two step learning (computed over $100\times 100$ points in the domain). 
    }
    \label{fig:mse_data_residual}
\end{figure}

We revisit the same linear quadratic toy example discussed in Section \ref{benefits}. The left-hand side of Fig.~\ref{fig:mse_training_methods} shows the MSE for three different training methods: i) training on value function data only, ii) training on data on the value function and the gradient, and iii) training on residual only. By minimizing over the residual only, the neural network converges to the incorrect solution (as was also shown in the toy example in Fig.~\ref{fig:lqr_loss_function}); by incorporating the gradient of the value function into the loss, a more accurate solution can be obtained. In the left-hand side of Fig.~\ref{fig:mse_data_residual} it is shown that two-step learning requires just 10 data points to obtain a solution that is of accuracy $10^{-6}$, while when learning only over the data, one would require 100 data points to achieve the same error. Comparing the left-hand sides of Fig~\ref{fig:lqr_loss_function} and ~\ref{fig:mse_data_residual} we also note the two-step learning can converge to a solution with a smaller MSE than the data-only training. Finally, the left-hand side of Table \ref{tb:data_vs_two_steps_results} shows the MSE for the data-only (using value function and gradient) and two-step learning. It is clear that when using more than 5 data points a more accurate solution is obtained using two-step learning. 

\begin{table}[hb]
\begin{center}
\caption{Mean squared error for 2D cases comparing the supervised learning on the dataset and the two steps learning.}\label{tb:data_vs_two_steps_results}
\begin{tabular}{l cc cc}
 & \multicolumn{2}{c}{2D linear} & \multicolumn{2}{c}{2D non linear} \\
\cmidrule(lr){2-3} \cmidrule(lr){4-5}
Point     & Data/grad   & 2 steps   & Data/grad   & 2 steps  \\ \hline
0 & \num{0.20288865} & \num{0.04377991} & \num{2.181217} & \num{1.176169} \\
5 & \num{0.00598487} &\num{0.02723277}  & \num{0.251505} &\num{1.179283} \\
10 &\num{0.00043885} & \num{0.00000109} &\num{0.04561} & \num{0.754814} \\
20 & \num{0.00004997} & \num{0.00000121} & \num{0.013527} & \num{0.000343} \\
50 & \num{0.00000280} & \num{0.00000153} & \num{0.00826} & \num{0.000206} \\
100 & \num{0.00000109} & \num{0.00000132} & \num{0.008112} & \num{0.000156}\\\hline
\end{tabular}
\end{center}
\end{table}

\subsection{A 2D nonlinear example}
\begin{figure}
    \centering
        \begin{tikzpicture}
            \begin{axis}[xlabel={Non linearity parameter $\epsilon$}, ylabel={Hamilton Jacobi residual}, width=0.9\linewidth, height= 0.6\linewidth,
            try min ticks=4,
            ylabel near ticks,
            xtick={0,1,2,3,4,5},
            xticklabels={$10^{-3}$, $10^{-2}$, $10^{-1}$, $1$, $10$, $100$}]
  
            \addplot[mark=x, mark size=4pt, color=black, thick]
            table[x index=0,y index=2,col sep=comma] {data/non_linear_discrepancy.dat};
            
            \end{axis}
        \end{tikzpicture}
       \caption{Influence of the intensity of the non linearity induced by $\epsilon$.}
       \label{fig:non_linear_discrepancy}
\end{figure}
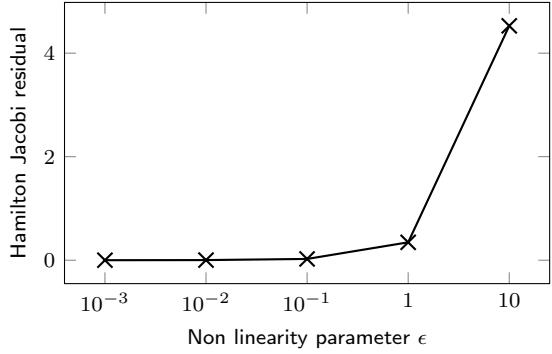

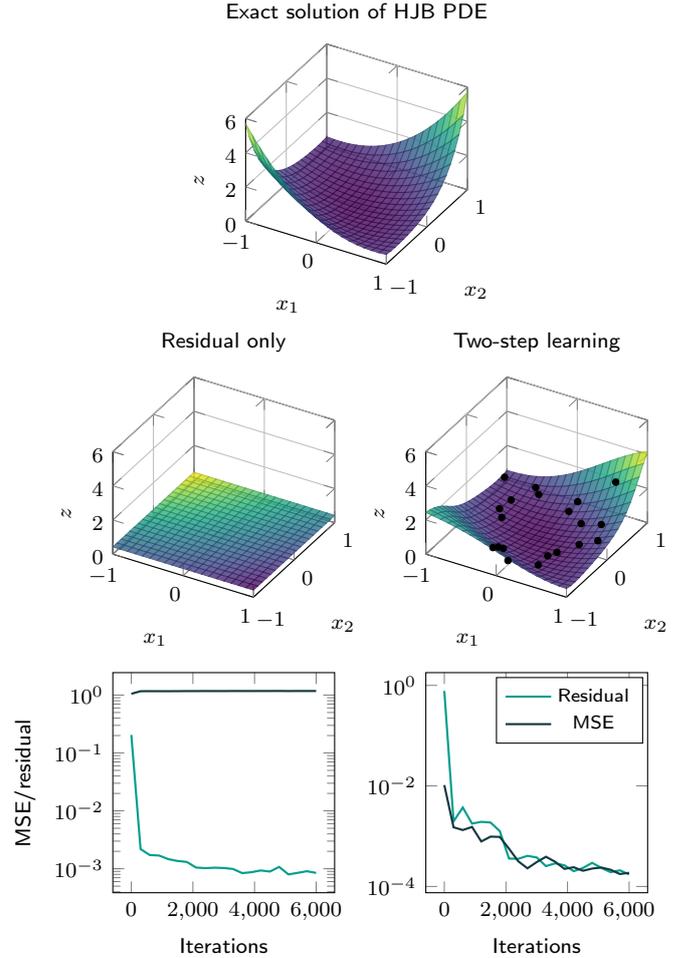
\begin{figure}
    \centering
\begin{tikzpicture}

\begin{groupplot}[group style={
        group name=my plots,
        group size=2 by 1,
        ylabels at=edge left,
    },
    width=4.5cm,
    height=4.5cm,
]

\nextgroupplot[
    title= Exact solution of HJB PDE,
    xlabel=$x_1$,
    ylabel=$x_2$,
    zlabel=$z$,
    view={30}{40},
    grid=both,
    mesh/ordering=x varies,
    mesh/rows=20,
    zmin=0,
    zmax=6.,
    ]
     \addplot3[
     opacity=0.8,
     colormap/viridis,
     surf,
     ] table[x index=0,y index=1, z index=6,col sep=comma] {data/function_all.dat};
     
\end{groupplot}
\end{tikzpicture}

\begin{tikzpicture}
\begin{groupplot}[
    group style={
        group name=my plots,
        group size=2 by 4,
        ylabels at=edge left,
        horizontal sep=1.2cm
    },
    width=4.5cm,
    height=4.5cm,
]

\nextgroupplot[
    title=Residual only,
    xlabel=$x_1$,
    ylabel=$x_2$,
    zlabel=$z$,
    view={30}{40},
    grid=both,
    mesh/ordering=x varies,
    mesh/rows=20,
    zmin=0,
    zmax=6.]
     \addplot3[
     opacity=0.8,
     colormap/viridis,
     surf,
     ] table[x index=0,y index=1, z index=2,col sep=comma] {data/function_all.dat};
         
\nextgroupplot[
    title= Two-step learning,
    xlabel=$x_1$,
    ylabel=$x_2$,
    zlabel=$z$,
    view={30}{40},
    grid=both,
    mesh/ordering=x varies,
    mesh/rows=20,
    zmin=0,
    zmax=6.]
     \addplot3[
     opacity=0.8,
     colormap/viridis,
     surf,
     ] table[x index=0,y index=1, z index=3,col sep=comma] {data/function_all.dat};
     
    \addplot3[only marks, mark size = 1.2] table[x index=0,y index=1, z index=2,col sep=comma] {data/dataset_coordinates.dat};

\nextgroupplot[ylabel=MSE/residual,
    ymode=log,
    xlabel=Iterations,]
    \addplot[color=graph_green, thick]
    table[x index=0,y index=1,col sep=comma] {data/non_linear_residual_solutions.dat};            
    \addplot[color=graph_petrol, thick]
    table[x index=0,y index=2,col sep=comma] {data/non_linear_residual_solutions.dat};  
            
\nextgroupplot[
    ymode=log,
    xlabel=Iterations,]
    \addplot[color=graph_green, thick]
    table[x index=0,y index=3,col sep=comma] {data/non_linear_residual_solutions.dat};            
    \addplot[color=graph_petrol, thick]
    table[x index=0,y index=4,col sep=comma] {data/non_linear_residual_solutions.dat};

    \addlegendentry{Residual}
    \addlegendentry{MSE}
\end{groupplot}
\end{tikzpicture}
    \caption{Training the residual from a random initialisation (left) and starting from a pre-training on a dataset (right). From top to bottom: the initialisation, the solution after training, the true solution and the value of the residual loss in \eqref{eq:residual_loss} and the MSE. 
    }
    \label{fig:invalid_vs_valid_residual}
\end{figure}

We next consider the dynamics given by
\[\dot \bx=\begin{bmatrix} 
	0 & 1  \\
	\epsilon x_1^2 & 0
\end{bmatrix}\bx+\begin{bmatrix} 
0  \\
1
\end{bmatrix}u(t)\,,\]
with $R = Q = I$. Note that the coefficient $\epsilon$ controls the intensity of the nonlinear term. The solution of the SDRE in this case is given by:
\begin{align*}
    P=\begin{bmatrix}
    \sqrt{\epsilon^2 x_1^4 + 1}\sqrt{1 + 2P_{12}} & P_{12} \\
    P_{12}  & \sqrt{1 + 2P_{12}}
    \end{bmatrix}
\end{align*}
with $P_{12} = \epsilon x_1^2 + \sqrt{\epsilon^2 x_1^4 + 1}$. 
Remember that the exact SDRE solution is an approximation of the original HJB PDE. In Fig.~\ref{fig:non_linear_discrepancy} we show the residual generated by the exact SDRE solution when inserted into the HJB PDE as a function of $\epsilon$ - this shows that the more non-linear the problem is, the more relevant the consequent minimization of the residual will be. This also shows why it is not possible to combine both \eqref{eq:residual_loss} and \eqref{eq:data_loss} into a single objective: the SDRE loss term and the residual loss term lead to different solutions and minimizing both at the same time leads to a decrease in accuracy compared to the two-step learning approach. 
The right-hand side of Fig.~\ref{fig:mse_training_methods} shows that using no data points and training on the residual only leads to inaccurate solutions, while training on the dataset only does not result in an error lower than $10^{-2}$ due to the SDRE solution being only an approximate solution. Combining data-driven initialisation with the residual minimization results in the most accurate solution as shown in the right-hand side of Fig.~\ref{fig:mse_data_residual}: training on 20 points first, and then minimizing the residual leads to a solution with a $10^{-4}$ accuracy. This level of accuracy is not obtained by increasing the amount of data points used in the supervised learning step. This highlights the necessity of the residual minimization to improve the accuracy of the final solution. The right-hand side of Table~\ref{tb:data_vs_two_steps_results} validates these results. Finally, in Fig.~\ref{fig:invalid_vs_valid_residual} we present the results from training on the residual only (the left-hand side) and the two-step learning (the right-hand side). The residual-only training results in an invalid solution. Based on the value of the residual loss in \eqref{eq:residual_loss} (the bottom plot) one could conclude that correct convergence has been obtained (the residual loss is sufficiently small); however as shown by the MSE computed by comparing the exact solution with the neural network output the network has converged to the incorrect solution. The value of two-step learning is clear from the right-hand side: the residual and MSE are both sufficiently low so that the solution obtained is close to the true solution. 

\subsection{Cucker-Smale model}
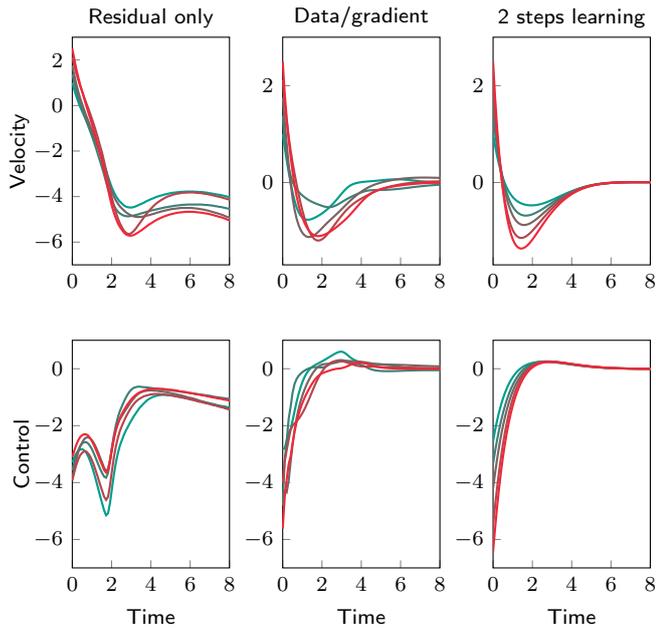
\begin{figure}
    \centering
\begin{tikzpicture}
\begin{groupplot}[
    group style={
        group name=my plots,
        group size=3 by 2,
        ylabels at=edge left,
        horizontal sep=0.7cm
    },
    footnotesize,
    width=3.65cm,
    height=4.6cm,
    tickpos=left,
    enlarge x limits=false,
    max space between ticks=1000pt,
    xmin=0,
    xmax=8,
]

\nextgroupplot[title={Residual only},ylabel=Velocity, ylabel shift = -5 pt, ymin=-7, ymax=3]
    \foreach [evaluate=\m as \redfrac using (\m-6)*100/(10-6)] \m in {6,...,10}{
        \edef\temp{%
            \noexpand
                \addplot[
                color=graph_red!\redfrac!graph_green,
                thick,
                ]
                table[x index=0,y index= \m] {data/pde_residual_invalid_trajectories.dat};
            }\temp
        }
\nextgroupplot[title={Data/gradient},ylabel shift = -5 pt, ymin=-1.7, ymax=3]
    \foreach [evaluate=\m as \redfrac using (\m-6)*100/(10-6)] \m in {6,...,10}{
        \edef\temp{%
            \noexpand
                \addplot[
                color=graph_red!\redfrac!graph_green,
                thick,
                ]
                table[x index=0,y index= \m] {data/pde_data_trajectories.dat};
            }\temp
        }
\nextgroupplot[title={2 steps learning},ymin=-1.7, ymax=3]
    \foreach [evaluate=\m as \redfrac using (\m-6)*100/(10-6)] \m in {6,...,10}{
        \edef\temp{%
            \noexpand
                \addplot[
                color=graph_red!\redfrac!graph_green,
                thick,
                ]
                table[x index=0,y index= \m] {data/pde_residual_trajectories.dat};
            }\temp
        }
\nextgroupplot[ylabel=Control, ylabel shift = -5 pt, xlabel=Time, ymin=-7, ymax=1]
    \foreach [evaluate=\m as \redfrac using (\m-1)*100/(5-1)] \m in {1,...,5}{
        \edef\temp{%
            \noexpand
                \addplot[
                color=graph_red!\redfrac!graph_green,
                thick,
                ]
                table[x index=0,y index= \m] {data/pde_residual_invalid_control.dat};
            }\temp
        }
\nextgroupplot[ylabel shift = -5 pt, xlabel=Time, ymin=-7, ymax=1]
    \foreach [evaluate=\m as \redfrac using (\m-1)*100/(5-1)] \m in {1,...,5}{
        \edef\temp{%
            \noexpand
                \addplot[
                color=graph_red!\redfrac!graph_green,
                thick,
                ]
                table[x index=0,y index= \m] {data/pde_data_control.dat};
            }\temp
        }
\nextgroupplot[xlabel=Time, ymin=-7, ymax=1]
    \foreach [evaluate=\m as \redfrac using (\m-1)*100/(5-1)] \m in {1,...,5}{
        \edef\temp{%
            \noexpand
                \addplot[
                color=graph_red!\redfrac!graph_green,
                thick,
                ]
                table[x index=0,y index= \m] {data/pde_residual_control.dat};
            }\temp
        }
\end{groupplot}
\end{tikzpicture}
    \caption{Trajectories and velocities of a group of particles under the control found by the neural network after (left) dataset-only training on 500 points and (right) pre-training on 500 points \& minimising the residual.
    }
    \label{fig:cuckersmale}
\end{figure}

We conclude with a high-dimensional nonlinear test consisting of agent-based consensus dynamics, known as the Cucker-Smale model. We consider 5 agents having states $\mathbf{x}_i=(y_i,v_i)\in\mathbb{R}^2$ in $[-3,3]^{10}\in\mathbb{R}^{5}\times\mathbb{R}^{5}$, with dynamics given in semilinear form by
\begin{align}
    \mathbf{\dot{x}}=\begin{bmatrix}
    O_5 & I_5\\
    O_5& A(\mathbf{y})
    \end{bmatrix} \mathbf{x} + \begin{bmatrix}
    O_5\\
    I_5
    \end{bmatrix}\bu\,,
\end{align}
where $O_5$ is an $5\times 5$ matrix of zeros and $I_5$ is the $5\times 5$ identity matrix and where
\begin{align}
A(\mathbf{y})_{ij}=\begin{cases}
	-\frac{1}{5}\sum_{k=1}^5 \frac{1}{1+||y_i-y_k||^2} & \text{ for } i=j\,,\\
	\frac{1}{5}\frac{1}{1+||y_i-y_j||^2} & \text{ otherwise}\,.
\end{cases}		
\end{align}
We set $Q = \frac{1}{5}I_{10}$ and $R=I_5$. Note that the resulting HJB PDE is 10-dimensional. In Fig.~\ref{fig:cuckersmale} we present the results for the position, velocity and control for i) residual-only training, ii) data and gradient-only training and iii) two-step learning. As before, it can be seen that training on just the residual  does not lead to the correct control. However, by first training on a dataset of 500 points, and then minimizing the residual, we are able to obtain a smooth stabilizing control. We remark that for this particular problem, 500 points is not enough for the supervised learning task, as in this data-only training the obtained control is not sufficiently smooth and many more points would be required to get a satisfying control that generalises well on unseen data. This highlights the advantage of our combined approach. 

\subsubsection{Conclusion and discussion}
We have shown how a two-step learning approach can result in efficient and accurate solutions to the HJB PDE. In our setting, the output of the neural network was set to directly be the value function $V(\mathbf{x})$. However, as in the SDRE derivation, it is possible to assume that $V(\bx)=\frac{1}{2}\mathbf{x}^{\top}P(\mathbf{x})\mathbf{x}$ and adapt the neural network architecture accordingly. Depending on the problem this can result in an efficient solver; we will address the choice of architecture for HJB PDE solvers in future work. Our results show that in a setting where the residual is sufficiently low for multiple solutions, the neural network converges to the one closest to initialization, which can be analyzed in the context of implicit bias \cite{chizat2020implicit}. 



\bibliography{ifacconf}             
                                                   








\end{document}